\documentclass[12pt]{article}
\usepackage{amsmath}
\usepackage{amsfonts}
\usepackage{a4wide}
\usepackage{amsthm}

\usepackage{amssymb}
\usepackage{latexsym}
\usepackage{euscript}

\newcommand{\Q}{\mathbf{Q}}
\newcommand{\Qbar}{\overline{\Q}}
\newcommand{\Qpbar}{\Qbar_p}
\newcommand{\Z}{\mathbf{Z}}

\title{A counterexample to the Gouv\^ea--Mazur conjecture.}
\author{Kevin Buzzard \and Frank Calegari\footnote{Supported in part by the American Institute of Mathematics.}}

\begin{document}
\maketitle

\begin{abstract}
Gouv\^ea and Mazur made a  precise conjecture about slopes of
modular forms. Weaker versions of this conjecture were
established by
Coleman and Wan. In this note, we exhibit examples contradicting
the full conjecture as it currently stands. 
\end{abstract}

\newtheorem{theorem}{Theorem}
\newtheorem{conjecture}[theorem]{Conjecture}
\theoremstyle{definition}
\newtheorem*{remark}{Remark}

Let~$p$ be a prime number, and let~$N$ be a positive integer coprime to~$p$.
For an integer $k$, let $f_k\in\Z[X]$ denote the characteristic polynomial
of the Hecke operator $U_p$ on the space of cusp forms of level $\Gamma_0(Np)$
and weight~$k$. Normalise the multiplicative valuation
$v:\Qpbar^\times\to\Q$ so that $v(p)=1$. If $\alpha\in\Q$ then
let $d(k,\alpha)$ denote the number of roots of $f_k$ in $\Qpbar$
which have valuation equal to $\alpha$.

\begin{conjecture}[Gouv\^ea--Mazur,\cite{gouvea-mazur}] If $k_1,k_2\in\Z$ are
both at least $2\alpha+2$ and $k_1\equiv k_2$~mod~$p^n(p-1)$ for
some integer $n\geq\alpha$, then $d(k_1,\alpha)=d(k_2,\alpha)$.
\end{conjecture}

This important conjecture was, as far as we know, one of the first precise
formulations of the idea that Hida's theory of ordinary families
might be generalised to the non-ordinary case. The conjecture
can be thought of as saying that non-ordinary families of cuspidal
eigenforms exist, and Gouv\^ea and Mazur went on to conjecture
that these conjectural families should in fact be analytic in the weight
variable. Many of the conjectures that Gouv\^ea and Mazur made
in this area became theorems soon afterwards, as a result of deep
ideas of Coleman (\cite{coleman:inventiones}), but the conjecture above,
which we shall refer to as the Gouv\^ea--Mazur conjecture, remained open.
Coleman's
ideas, extended by Wan in~\cite{wan}, could only prove that if $k_1$ and $k_2$
were congruent modulo $p^N(p-1)$ for some integer $N=O(\alpha^2)$
then $d(k_1,\alpha)=d(k_2,\alpha)$. For a while the situation
was almost paradoxical, because computations of Gouv\^ea, Stein,
and the first author all seemed to indicate that in fact the Gouv\^ea--Mazur
conjecture was much too weak. See for example~\cite{buzzard:conjectures}
which introduces the notion of a prime $p$ being $\Gamma_0(N)$-regular
and makes some much more precise conjectures when this condition holds.

If $N=1$ then the smallest prime which is not $\Gamma_0(N)$-regular
is $p=59$, and the authors decided to attempt to study this case
in detail. In particular, they wanted to study the point of the
eigencurve corresponding to the slope~1 form of level~59
and weight~16. To their surprise, they discovered

\begin{theorem} Set $N=1$ and $p=59$. Then there exists $\alpha\in\Q$
with $0\leq\alpha\leq 1$ such that $d(16,\alpha)\not=d(3438,\alpha)$.
In particular, because $3438=16+58 \cdot 59$, the Gouv\^ea--Mazur conjecture
is false.
\end{theorem}
\begin{remark} We know $d(16,1)=1$, and it is probably true that $d(3438,1)=2$,
but it is, at the present time,
difficult to compute characteristic polynomials of Hecke
operators at such high weight.
\end{remark}
\begin{proof} It is well-known that the space of level~1 cusp forms
of weight~16 is 1-dimensional, and that the 59-adic valuation
of the eigenvalue of $T_{59}$ is~1 on this space. Because any
newforms of level~59 and weight~16 will have slope $(16-2)/2=7$,
we see that $d(16,1)=1$ and $d(16,\alpha)=0$ for $0\leq\alpha<1$.

Now using the Eichler--Selberg trace formula (see for
example
Theorem~2 on p48 of~\cite{lang:modularforms}, and~\cite{zagier:correction}),
it is possible
to compute the trace of $T_{59}$ and $T_{59^2}$ on the space of cusp forms
of level~1 and weight~3438. These can be computed as real numbers to
very large precision in about 30 seconds using {\tt pari-GP}
 on a 300MHz Pentium III
(the trace of $T_{59^2}$ is about $6.79926\ldots\times 10^{6086}$),
and these real numbers turn out to be very close to integers, so one
can round to get the exact result. In fact, for the argument below, we
only need to know the traces modulo $59^3$ and these can be computed in 
less than a tenth of a second using a modular implementation of
the trace formula.

Let the characteristic polynomial of $T_{59}$ on the space
of cusp forms of level~1 be $X^{286}+\sum_{i=1}^{286}a_iX^{286-i}$.
From the above calculations one has enough data to compute
$a_1$ and $a_2$ exactly (although we only need them modulo $59^3$).
One now checks that the 59-adic valuations of $a_1$ and $a_2$
are~1 and~2 respectively. There are now two cases to consider:
either $v(a_i)\geq i$ for all~$i$, in which case $d(3438,1)\geq2$,
or there exists $i$ with $v(a_i)<i$, in which case $d(3438,\alpha)>0$
for some $\alpha$ with $0\leq\alpha<1$. In either case the conjecture
is violated.
\end{proof}
\begin{remark} The authors strongly suspect that $v(a_i)\geq i$ for
all~$i$, and perhaps a more intense computation and some estimates
on the $a_i$ might prove this.
\end{remark}

This counterexample comes from a point on an eigencurve with
integral slope, whose associated mod~$p$ representation is irreducible
even when restricted to a decomposition group at~$p$. In fact
a search for other similar points led to other counterexamples:

\begin{theorem} Set $p=5$. Then
the set of $N$ with $1\leq N\leq 83$ and such that $5\nmid N$
and $d(6,1)>0$ is $\{14,28,34,37,38,42,53,56,68,69,71,74,76,83\}$,
and for each $N$ in this set we have $d(26,1)=2d(6,1)>d(6,1)$,
so every such~$N$ gives a counterexample to the Gouv\^ea--Mazur
conjecture.
\end{theorem}
\begin{proof} One verifies this statement easily using William
Stein's programs for modular forms, written for the {\tt MAGMA}
computer algebra package~\cite{magma}.
\end{proof}
\begin{remark} As opposed to the $p=59$ counterexample above,
which relies on ``custom computations'' done by the authors,
these latter counterexamples rely only on standard programs of Stein
which have existed for many years and are likely to be
bug-free.
\end{remark}

The
counterexamples exhibited above all pertain to forms
of positive slope and weight at most $p+1$, and hence to
forms whose associated mod~$p$ representation is irreducible
when restricted to a decomposition group at~$p$. After our
counterexamples were discovered, Graham Herrick pointed out to us
that if one extends the domain of the conjecture to forms of
level $\Gamma_1(N)\cap\Gamma_0(p)$ then one can find counterexamples
coming from forms whose associated mod~$p$ representation is reducible,
even globally. For example if $p=5$ then there are two slope~1 forms
of level $\Gamma_1(4)$ at weight~7, and four slope~1 forms of this
level at weight~27. The mod~5 representations associated to all
of these forms are reducible. Note that once again there are twice as many
slope~1 forms in the higher weight than are predicted by the Gouv\^ea-Mazur
conjecture.

Finally, we remark that Breuil has made local conjectures
(Conjecture~1.5 in~\cite{breuilII}) which imply that if~$f$
is an eigenform of level~$N$ and weight~$k$ with $2\leq k\leq 2p$,
and if the slope of~$f$ is strictly between~0 and~1, then the
mod~$p$ representation associated to~$f$ is irreducible
when restricted to a decomposition group at~$p$. Breuil has
asked (personal communication with K.B.) whether one should
expect such a result when $k>2p$. However, the mod~5 representations
associated to the four forms of level $\Gamma_1(4)$ and weight~11
are all reducible, and all of these forms have slope~$1/2$. The
generalisation of Breuil's conjectures to the case $k>2p$ remains
an interesting open problem, which might shed more light on the
results of this paper.

{\tt k.buzzard@ic.ac.uk} \\
{\tt fcale@math.harvard.edu}
\end{document}